\documentclass[11pt, 
]{amsart}
\usepackage{ 
amssymb 
} 
\newcommand{\no}[1]{#1} 

\renewcommand{\no}[1]{} 
 
\no{\usepackage{times}\usepackage[ 
subscriptcorrection, slantedGreek,  
nofontinfo]{mtpro} 
\renewcommand{\Delta}{\upDelta} 
} 
 
\date{\today}

\newtheorem{theorem}{Theorem}[section]

\theoremstyle{remark} 
\newtheorem{remark}{Remark}[section]

\newcommand{\R}{{\bf R}} 
 
\renewcommand{\r}[1]{(\ref{#1})} 
\newcommand{\PDO}{$\Psi$DO} 
\newcommand{\be}[1]{\begin{equation}\label{#1}} 
\newcommand{\ee}{\end{equation}} 

\renewcommand{\d}{\mathrm{d}}

\newcommand{\bo}{{\partial \Omega}}

\numberwithin{equation}{section}

 \usepackage{graphicx}
 
\title[Modulated Luminescent Tomography ]{Modulated Luminescent Tomography} 
 
\author[P. Stefanov]{Plamen Stefanov} 
\address{Department of Mathematics, Purdue University, West Lafayette, IN 47907} 
\thanks{First author partly supported by the NSF  Grant DMS-0800428} 
\author[Wenxiang Cong]{Wenxiang Cong} 
\address{Department of Biomedical Engineering, Rensselaer Polytechnic Institute, Troy, NY 12180}

\author[Ge Wang]{Ge Wang} 
\address{Department of Biomedical Engineering, Rensselaer Polytechnic Institute, Troy, NY 12180}
\thanks{Second and third author partly supported by the NIH Grant NIH/NHLBI HL098912, and the NSF Grants CMMI-1229405 and DMR-0955908}

\begin{document} 
 \begin{abstract} We propose and analyze a mathematical model of Modulated Luminescent Tomography. We show that when single X-rays or  focused X-rays are used as an excitation, the problem is similar  to the inversion of  weighted X-ray transforms. In particular, we give an explicit inversion in the case of Dual Cone X-ray excitation.
 \end{abstract}

 \maketitle 
\section{Introduction: A brief description of the method and the model}
\subsection{ A brief description of the modalities} \ 
The purpose of this work is to build and analyze the mathematical model of several medical multiwave imaging techniques that we call Modulated Luminescent Tomography. Certain phosphorescent contrast agents (for example, nanophosphors) are delivered to the cells and then illuminated with    X-rays. The illuminated particles emit photons which scatter through the tissue and are detected outside the body, see, e.g., \cite{Liu2013}. The goal is to recover the concentration of the contrast agents. The resolution is expected to come from the X-ray excitation. In this respect, those modalities behave differently than Photoacoustic and Thermoacoustic Tomography, where the excitation is highly diffusive but the emitted ultrasound signal allows for high resolution, see, e.g., \cite{Wang2012}.

Two such techniques, X-ray luminescence computed tomography (XLCT),  \cite{Pratx10} and X-ray micro-modulated luminescence tomography (XMLT) \cite{Cong14, St-double_cone_eng}, described below, have been proposed recently. 

Microscopy is the principal observational tool that has made fundamental contributions to our understanding of biological systems and engineered tissues. Popular microscopy techniques use visible light and electrons. Sample preparation and imaging with these techniques are relatively simple, being good for in situ or in vivo studies of cultured cell/tissue samples. Inherently, the resolution of optical microscopy is limited by diffraction. With additional sample preparation, stochastic information and innovative interference techniques ~100nm resolution is achievable. Three-dimensional images can be obtained with optical sectioning. Ultimately, multiple scattering prevents these techniques from imaging thick samples. Photoacoustic tomography permits scalable resolution with a depth of up to ~7cm with a depth-to-resolution ratio ~200. Photoacoustic microscopy aims at millimeter depth and micron-scale resolution based on absorption contrast \cite{Wang2013}, which can be used to characterize the structure of a scaffold but does not provide the sensitivity of fluorescence imaging. All these methods are bounded by 1mm imaging depth.   

In the in vivo imaging field, fluorescence molecular tomography (FMT) and bioluminescence tomography (BLT) are capable of visualizing biological processes at molecular and cellular levels deep inside living tissue \cite{Ntziachristos2005}. Such optical molecular imaging tools have major advantages in terms of high sensitivity, biological specificity, non-invasiveness, and cost-effectiveness, which are widely used in preclinical investigations and a limited number of clinical applications \cite{Wang2006}. However, resolution and stability of FMT and BLT still remain unsatisfactory. Recently, X-ray luminescence computed tomography (XLCT) was developed \cite{Pratx10} with nanophosphors (NPs) as imaging probes. Excited with a pencil beam of X-rays, NPs simultaneously gives rise to luminescence emission whose photons can be efficiently collected with a CCD camera. This mechanism allows tomographic reconstruction of a NP distribution, similar to the case of X-ray CT. X-ray excitation not only confines luminescence sources within the active beam but also eliminates auto-fluorescence related artifacts in FMT images.  However, due to several technical and physical limitations, the lower bound of XLCT resolution is about one millimeter and seems infeasible to break through.

Recently,  a unique imaging approach called “X-ray micro-modulated luminescence tomography (XMLT)” approach was proposed \cite{Cong14, St-double_cone_eng}, which combines X-ray focusing, nanophosphor excitation, optical sensing, and image reconstruction in a synergistic fashion, and promises significant imaging performance gains in both microscopic and preclinical studies. The XMLT approach is
an exogenous reporter-based imaging system in which spatial resolution is determined by the excitation of the luminescence from nanophosphors with a micro-focused X-ray beam such as a zone-plane or a micro-focus X-ray source coupled with a polycapillary lens. High in vivo spatial resolution can be achieved due to the short wavelength and small spot size of the X-rays.  When X-rays are focused in this way, double cones of X-rays are formed with their shared vertex point inside a sample or a subject.

Our main results state that we can use boundary averaged measurements only and the problem is reduced to an inversion of a weighted linear transform determined by the excitation. For example, if the excitation consist of single X-rays, in Section~\ref{sec5}, we get a weighted X-ray transform with a weight depending on the diffusion and the absorption coefficients $D$ and $\mu_a$ of the medium with respect to the emitted photons, see \r{M1} below.  
In the case of double cone excitation (XMLT) with a constant aperture, or in XLCT, we have an explicit and stable inversion under an explicit  if and only if condition \r{stab_cond}, see Theorem~\ref{thm4.1} and Section~\ref{sec4}. We show that the local problem in a region of interest behaves similarly to the corresponding X-ray problems, and we have explicit microlocal inversions of the visible singularities.

\subsection{The model}   

Let $\Omega\subset \R^n$ be a smooth domain. 
We illuminate the medium in different ways. Each illumination excites the particles at a rate $I_\alpha(x)$, which can be a distribution, where $\alpha$ belongs to some index set $\mathcal{A}$, discrete or not. In the X-ray case, $I_\alpha(x)$ is just a superposition of delta functions along straight lines; which also allows for  attenuated X-rays by introducing an appropriate weight. The structure of $I_\alpha$, XLCT or XMLT or something else, see section~\ref{sec4} and section~\ref{sec5}, is not important at the moment.  
The modulated source then is
\be{S}
s_\alpha(x) = I_\alpha(x)f(x),
\ee
where $f$ is the concentration of the active particles, which we want to reconstruct.   We assume that the light from the source propagates according to the diffusion model \cite{Welch-book}
\be{M1}
\left(-\nabla\cdot D\nabla +\mu_a\right) u_\alpha= s_\alpha (x) \quad \text{in $\Omega$}. 
\ee
Here, $\mu_a(x)\ge0$ is the absorption, $D(x)>0$ is the diffusion coefficient,  
and $u_\alpha$ is the photon density. The coefficients $\mu_a$ and $D$ are given. In the diffusion approximation regime, $D(x)= [3(\mu_a+\mu'_s)]^{-1}>0$, where   $\mu_s'=(1-g)\mu_s$ is the reduced scattering coefficient, $\mu_s$ is the scattering coefficient and $g$ is the scattering anisotropy. 

The boundary conditions are of semi-transparent type due to the different index of refraction of the tissue and the air around it. They are of Robin type and have the form
\be{M2}
u_\alpha+2AD\partial_\nu u_\alpha \big|_{\bo}=0,
\ee
where $A>0$ is a given coefficient, well approximated by $A=(1+R)/(1-R)$ with $R$ closely approximated by $R= -1.4399m^{-2} + 0.7099m^{-1}+ 0.6681+0.063m$ \cite{Schweiger95}, where $m$ is the  refractive index   of the tissue.  
Let $G$ be the solution operator (the Green's function) of \r{M1}, \r{M2}, i.e., $u_\alpha = Gs_\alpha = GI_\alpha f$. 

What we measure is the outgoing photon density
\be{M3}
Q_\alpha f= -D\partial_\nu u_\alpha \big|_{\bo} =\frac1{2A} u_\alpha |_{\bo} =\frac1{2A} (GI_\alpha f)|_{\bo}.
\ee

The data therefore is
\[
\{Q_\alpha f\}_{\alpha\in\mathcal{A}},
\]
and we want to reconstruct $f$. 
For each illumination choice, $\alpha$, the data $Q_\alpha$ is a function of $n-1$ variables. When $\alpha$ runs over a continuous space with dimension $m$, we get $m+n-1$ variables. This could make the problem formally overdetermined. On the other hand, because of the diffusion nature of the data for a fixed $\alpha$, we cannot have much resolution hidden in that $Q_\alpha$. The only way to get resolution (i.e., stability) is to have a well chosen set of illuminations $I_\alpha$,  $\alpha\in\mathcal{A}$ which have enough singularities. We show below that we do not actually need to know or measure $Q_\alpha$ pointwise; some average over $\bo$ is enough for a stable recovery. This removes $n-1$ variables from the data. 

\section{The first phase: recovery of an averaged intensity.} 
The operator $L:= -\nabla\cdot D\nabla +\mu_a$ is symmetric and positive on smooth functions satisfying the the Robin boundary condition \r{M2} because
\be{L}
\begin{split}
\int_\Omega (Lu)\bar u\,\d x &= \int_\Omega\left(  |\nabla u|^2+\mu_a|u|^2 \right)\d x - \int_\bo (\partial_\nu u)\bar u\,\d\sigma\\
 &=  \int_\Omega\left(  |\nabla u|^2+\mu_a|u|^2 \right)\d x +\int_\bo \frac{1}{2AD}  |\bar u|^2\,\d\sigma,
\end{split}
\ee
where $\d \sigma$ is the surface measure.  By 
Green's formula, for such $u$ and $v$,
\be{M4}
\begin{split}
\int_\Omega \left( vLu -uLv\right)\,\d x &= -\int_{\bo} \left(  v\partial_\nu u - u\partial_\nu v    \right)D\,\d \sigma\\
 &= -\int_{\bo} \big(  (v +2AD\partial_\nu v) \partial_\nu u - (u  +2AD\partial_\nu u )\partial_\nu v    \big)D\,\d \sigma    =0. 
\end{split}
\ee
Therefore, $L$ with the Robin boundary conditions if symmetric and positive by \r{L}. 
It is well known that $L$ extends to a self-adjoint operator on $L^2(\Omega)$  with a compact resolvent. In particular, $0$ is in the resolvent set by \r{L}. This shows that the Green's function $G$ is well-defined, and self-adjoint on $L^2(\Omega)$. Moreover, the boundary-value problem \r{M6} below is well posed. 

If $Lv_\alpha=0$, and $u_\alpha$ solves  \r{M1}, \r{M2},   we get (we drop the index $\alpha$ in this formula)
\be{M5}
\begin{split}
\int_\Omega   vs\,\d x &= -\int_{\bo} \left(  v\partial_\nu u - u\partial_\nu v     \right)D\,\d \sigma= -\int_\bo \left(   v +2AD \partial_\nu v    \right)D\partial_\nu u \,\d \sigma \\
&= \int_{\bo} \left(   v +2AD \partial_\nu v    \right)Qf \,\d \sigma. 
\end{split}
\ee
This suggests the following. Choose $h_\alpha$ and let $v_\alpha:=Vh_\alpha$ solve
\be{M6}
\begin{split}
\left(-\nabla\cdot D\nabla +\mu_a\right) v_\alpha&= 0 \quad \text{in $\Omega$},\\
 v_\alpha+2AD\partial_\nu v_\alpha\big|_{\bo}&=h_\alpha. 
\end{split}
\ee
 Then, by \r{M5},
\be{M7}
\int_\Omega v_\alpha  I_\alpha f  \,\d x = \int_{\bo}h_\alpha Q_\alpha f\,\d \sigma.
\ee
In other words, choose $h_\alpha$ somehow; then we can recover $\int_\Omega\! (Vh_\alpha)I_\alpha f\,\d x$ from the data. This requires us to solve \r{M6} first but this can be done numerically in a very efficient way. If $D$ and $\mu_a$ are constants, and $\Omega$ is a circle or a rectangle, it can be done explicitly, as we show below. This gives us a way to compute an averaged value of $I_\alpha f$, with a weight $v_\alpha$ depending on the choice of $h_\alpha $. If   $\mu_a=0$, one can take $h_\alpha=1$ and then $v_\alpha =1$. 
 We would want to have positive solutions $v>0$, and we do not want $v\ll1$ because we will divide by it eventually. Positivity can be guaranteed by the maximum principle if $h>0$, since  $0<v<\max_{\bo}h$. 
The lower bound of $u$ could be very small if the attenuation $\mu_a$ is large, and that would create stability problems; but this is natural to expect. 

Note that we do not need to solve the boundary value problem \r{M6} to find $v$. We could just take any solution $v$ of $(-\nabla\cdot D\nabla +\mu_a)v= 0$ and compute $h=  (v+2AD\partial_\nu v)\big|_{\bo}$ next. If $D$ and $\mu_0$ are constants, we can generate many explicit solutions, regardless of the shape of the domain.

Finally, we want to emphasize that we cannot recover $s_\alpha= I_\alpha f$ from \r{M7} known for all $h_\alpha$; for example,  for any $\phi\in C_0^\infty(\Omega)$, adding $L\phi$ to $I_\alpha f$ would not change the left-hand side of \r{M7}.

We consider the constant coefficient cases below and continue with the general case in the next section. 

\subsection{Partial cases: $D$ and $\mu_a$ constants}
Assume now that $D>0$ and $\mu_a\ge0$ are constants. 

\subsection{$n=3$, $\Omega$ a ball} 
Let $n=3$. The Green's function of $L=-D\Delta+\mu_a$ is then given by 
\be{MG}
G_k(x,y) = \frac{e^{-k|x-y|}}{4\pi D|x-y|},
\ee
with $k=\sqrt{\mu_a/D}$, i.e., $(-D\Delta_x+\mu_a) G(x,y)=\delta(x-y)$ 
Then $G_{-k}$ is a Green's function as well, and the difference is in the kernel of $L$. Therefore,
\[
(-D\Delta_x+\mu_a)\frac{\sinh\left(k|x-y|\right)}{|x-y|}=0. 
\]
Therefore, we can chose
\be{Mv}
v(x) = \frac{\sinh\left(k|x-y|\right)}{|x-y|}
\ee
with $y$ arbitrary but fixed. Then compute $h$ from \r{M6}.

If $\Omega$ is the ball $\Omega= \{x|\;|x|\le a\}$, then we can choose $y=0$ to get
\be{M8a}
v(x) = \frac{\sinh{(k|x|)}}{|x|}, \quad 
h=  a^{-1}\sinh(ka)+     2AD\left(ka^{-1}\sinh(ka)- a^{-2}\sinh(ka)\right). 
\ee
Note that $h$ is just a constant on the boundary $|x|=a$. 
Then by \r{M7}, we can recover 
\[
\int_{|x|\le a}I_\alpha(x)f(x)\frac{\sinh{(k|x|)}}{|x|}\, \d x. 
\]
When $ka\gg1$, the function $v$ is very small near $0$ compared to its boundary values, which means poor recovery of $s_\alpha$ there (we do not actually recover $s_\alpha$ here, just a weighed integral of it). That can be expected though --- signals coming from the center would  be attenuated the most.

\subsection{$n=3$, $\Omega$ arbitrary} One can still use $v$ as in \r{M8a} (or, as in\r{Mv}, see the remark above); then $h$ needs to be calculated  by the second equation in \r{M6}. 

\subsection{$n=2$, $\Omega$ a disk} In this case, the fundamental solution involves the Bessel function $K_0$ and is $K_0(k|x-y|)$. Note that there is a weak (logarithmic) singularity at $x=y$. This is equivalent to \r{MG} in the 3D case.  One can then take
\[
v(x) = I_0(k|x|),
\]
where $I_0$ is a Bessel function again, and there is no singularity. This is equivalent to taking $\sinh(k|x|)/|x|$ in 3D. For $h$ on the sphere $|x|=a$ the boundary of the ball), we get
\[
h= I_0(ak)+ 2ADI_1(ak),
\] 
which is a constant again.

\section{Reducing the problem to the invertibility of a weighted version of the illumination transform}\label{sec3} 
As we showed above, we cannot recover the source $s_\alpha =I_\alpha f$ from \r{M7}. We are not trying to reconstruct $s_\alpha$ however. 
 We view \r{M7} as a linear integral transform 
\be{M8}
R : f \longmapsto Rf(\alpha) := \int_\Omega v_\alpha I_\alpha f\, \d x, \quad \alpha\in \mathcal{A}.
\ee
We can write equation \r{M7} then as
\be{Q1}
Rf(\alpha) = \int_\bo  h_\alpha Q_\alpha f\,\d\sigma.
\ee
The right-hand side is determined by the data. 
The goal then is to invert $R$. Recall that we have some freedom to choose $v_\alpha=Vh_\alpha$, and in particular, we can choose them independently of $\alpha$. Then $v$ just multiples $f$, and we need to invert the transform with kernel $I_\alpha(x)$; and then divide by $v$. 

We consider two main examples below: when $R$ is a dual cone transform modeling the XMLT; and when $R$ is a weighted X-ray transform, which models the XLCT.

\section{Double cone excitation (XMLT)}\label{sec4}
\subsection{Formulation} 
 We work in $\R^n$ but the interesting case is $\R^3$. 
The idea of the XMLT is to focus X-rays, thus forming a double cone, 
at each point $x$ in some region of interest ROI at $N$ different directions $\theta_j$, $j=1,\dots, N$. Then 
\[
\alpha=(x,j)\in \Omega\times \{1,2,\dots,N\}.
\]
We model the cones with their aperture functions $a_{x,j}(\theta)$, $\theta\in S^{n-1}$, in other words, $I_\alpha(x)$ is a superposition of X-rays with density $a_{x,j}$. 
 Then the corresponding  intensity $I_{x,\alpha}(y)$ is given by
\[
I_{x,j}(y) = \frac{  a_{x,j}\big(\frac{x-y}{|x-y|}\big)}{|x-y|^{n-1}}.
\]
For example, if we assume uniform angular density, then $a_{x,j}$ would be constant on the intersection of the interior of the cone with the unit sphere, and zero otherwise. We assume that the cones are double; then $a_{x,j}$ are even functions on the sphere. The simplest case is when $a_{x,j}$ is independent of $x$, and $a_{x_1,j_1}$ is obtained from $a_{x_2,j_2}$ by a translation in the $x$ variable and a rotation in the angular one. 

By \r{S}, the transform $R$ becomes
\be{X2a}
Rf(x,j) =\int_\Omega\frac{  a_{x,j}\big(\frac{x-y}{|x-y|}\big)}{|x-y|^{n-1}}v_{j,x}(y)f(y)\,\d y.
\ee
If $a_{x,j}$  are smooth functions, then $R$ is a \PDO\ of order $-1$ with principal symbol
\[
r_j(x,\xi) := 
\pi \int_{S^{n-1}}a_{x,j}(\theta)   v_{j,x}(x)  \delta (\xi\cdot\theta)\,\d\theta,
\]
see, e.g., \cite{FSU}. Here, $\delta$ is the Dirac delta ``function''; therefore, the integral is taken over the sphere of co-dimension two (the grand circle if $n=3$, two points if $n=2$) of the unit sphere intersected with the plane normal to $\xi$. 
Note that the symbol above is homogeneous of order $-1$. This shows that the ``visible'' set of singularities $\mathcal{U}\subset T^*\Omega\setminus 0 $ is where the system $\{r_j\}_{j=1}^N$ is elliptic, i.e.,
\be{U}
\mathcal{U} = \{  (x,\xi)\in T^*\Omega\setminus 0;\;  \exists j,   \theta\perp\xi\; \text{so that}\; a_{x,j}(\theta)\not=0   \}.
\ee
\subsection{Recovery of the singularities in a region of interest} If $\Omega_0\subset\Omega$ is an open subset (a region of interest), and if $ a_{x,j}$ are smooth, we can recover the singularities of $f$ there, if the condition for $\mathcal{U}$ in \r{U} is satisfied for any $x\in\Omega_0$, i.e., if 
\be{stab_cond}
T^*\Omega_0\setminus 0\subset\mathcal{U}.
\ee
 This can be done constructively as follows.
\begin{itemize}
\item For any $j$ and $x\in \Omega_0$,   choose a smooth function $0<h_{j,x}$ on $\bo$ and compute $v_{j,x}=Vh_{j,x}$ by solving \r{M6} with $h_\alpha =h_{j,x}$.  This determines the transform $R$. 
\item Compute the right-hand side of \r{Q1} 
 for $x\in\Omega_0$. 
\item In \r{Q1}, apply a left parametrix for the matrix operator $f\mapsto Rf(j,\cdot)$ in $\Omega_0$. 
\end{itemize}
One way to construct a parametrix of order $1$ is to construct a \PDO\ with the following principal symbol
\[
q_j := \Big (\sum_kr_k^2\Big)^{-1}   r_j,
\]
and then $\sum_j  q_j(x,D)Rf(j,x)=f+Kf$, where $K$ is of order $-1$. If $a_{x,j}\ge0$, one can use instead
\[
 q := \Big (\sum_k r_k\Big)^{-1}
\]
and then $ q(x,D)\sum_j Rf(j,x)=f+Kf$ with $K$ as above. 

Note that we can simply choose $h_{j,x}$ independent of ${j,x}$.

The aperture functions $a_{x,j}$ only need to be smooth of finite order to recover the leading singularities, like jumps across surfaces. In practical implementations, this condition can be satisfied approximately by introducing  vignetting near the edges of the cones.  Even without applying the operator $Q$, the $\sum_j Rf(j,x)$ would recover all singularities in the correct places but with varying amplitudes. 

When the stability condition \r{stab_cond} does not hold, we can still reconstruct the visible singularities in $\mathcal{U}$ stably  by a microlocal inversion, if $a_{x,j}$ are smooth but  not the ones in the complement of $\bar{\mathcal{U}}$, see, e.g., \cite{SU-JFA09}. 

\subsection{Explicit global inversion when the apertures are translation-independent} 
Assume now that the aperture functions are non-negative and independent of $x$, i.e., they are given by $a_j(\theta)$.  The stability condition \r{stab_cond} then takes the following form, see \r{U} :
\be{U2}
\text{Any plane through the origin intersects at least one of the cones $\{\theta;\; a_j(\theta)>0\}$.}
\ee
Assume also that $h_{j,x}$ are independent of $j$ and $x$.  
The construction above then simplifies as follows. We have
\[
r_j(x,\xi) :=   r^0_j(\xi)v(x), \quad r^0_j(\xi):= \pi 
\int_{S^{n-1}}a_j(\theta)    \delta (\xi\cdot\theta)\,\d\theta,
\]
and $r_j(\xi)v(y)$ is actually the full (not just the principal symbol) amplitude of $f\mapsto Rf(j,\cdot)$. Symbols independent of $x$ are  Fourier multipliers, i.e., $r_j(D)= \mathcal{F}^{-1}r_j\mathcal{F}$, where $\mathcal{F}$ is the Fourier transform. 
Then $R_j= r^0_j(D)v$ with $v$ regarded as a multiplication, where we changed the notation for $R$ since $\alpha=j$ now. In integral form,
\[
R_jf(x) = \int_\Omega\frac{  a_j\big(\frac{x-y}{|x-y|}\big)}{|x-y|^{n-1}}v(y)f(y)\,\d y. 
\]
Note that $r_j^0(D)$ makes sense even if $a_j$ is $L^\infty$ only.

Equation \r{Q1} takes the form
\[
R_jf (x)= \int_\bo h Q_{j,x}\,\d\sigma. 
\]
When the stability condition \r{U2} holds,   this equation can be solved explicitly as follows
\be{f}
f = \frac1v \Big(\sum_j r_j(D)  \Big)^{-1}\sum_j \int_\bo h Q_{j,x}\,\d\sigma. 
\ee
We summarize this in the following:

\begin{theorem}\label{thm4.1}
Let the aperture functions $0  \le  a_j(\theta)\in L^\infty$ be independent on the point $x$ where we focus.  Choose $0<h\in L^\infty(\bo)$ and let $v=Vh$ be the solution of \r{M6} corresponding to that $h$. If the stability condition \r{U2} holds, then $f$ can be explicitly computed by \r{f}.  

Moreover,
\[
\|f\|_{L^2(\Omega)} \le C\frac{\sup_{\bo}h}{\inf_\Omega v}\Big\| \int_\bo \sum_j Q_{j,x}\,\d\sigma  \Big\|_{H^1(\Omega_x)}.
\]
\end{theorem}

The subscript $x$ in $\Omega_x$ indicates that the $H^1$ norm is taken w.r.t.\ to the variable $x$. 

\begin{remark}\label{rem1}
One way to satisfy the stability condition \r{U2} is to make sure that union of the open sets $\{a_j>0\}$  where the apertures do not vanish covers the equator of the unit sphere (or any fixed in advance grand circle). This is the situation in the numerical example below. If the apertures $\{a_j>0\}$ are too small, then there will be a large variation between the minimum and the maximum of the symbol $q$ below, which may lead to a mild instability; the most stable singularities would be the vertical ones. 
\end{remark}

\section{X-ray excitation (XLCT)} \label{sec5}
In the X-ray excitation case, we send individual  X-rays through the medium. If the direction is $\theta\in S^{n-1}$, and we parameterize lines in this direction by initial points $z$ on $ \theta^\perp =\{x;\; x\cdot\theta=0\}$, we have
\[
\alpha = (\theta,z), \quad \theta\in S^{n-1}, \; z\in \theta^\perp. 
\]
In other words, we identify the directed lines in $\R^n$ with such $\alpha$'s. Then $I_\alpha$ is just a delta unction on the line parameterized by $\alpha$. Therefore, $R$ is just the weighted X-ray transform with weight $h_{\theta,z}$. Assume now that $h$ is chosen independently of the line parameterized by $(\theta,z)$. Then
\[
Rf= R_0vf, \quad R_0g(\theta,z) = \int g(z+t\theta)\,\d t.
\]
Equation \r{Q1} takes the form
\[
(R_0vf)(\theta,z)= \int_\bo hQ_{\theta,z}\,\d\sigma. 
\]
In other words, each measurement gives us an integral of $vf$ over a single line. 
The function $f$ can then be found by inverting the X-ray transform $R_0$ and then dividing by $v$. The problem with a subset of lines and microlocal recovery in a ROI is the same as with the X-ray transform. 
Note that in practice, diffraction and other engineering challenges limit the ability to concentrate the radiation too closely along a single line thus limiting the resolution.

\section{Numerical simulations}

A spherical phantom of radius 10 mm was employed for the numerical simulation. The phantom is assigned with optical parameters: absorption coefficient $\mu_a=0.05$ mm$^{-1}$, scattering coefficient $\mu_s=15.$0 mm$^{-1}$, anisotropic coefficient $ g=0.9$, and relative refractive index of $1.37$. The phantom was discretized into $65,775$ tetrahedral elements and $12,044$ nodes. A total of $2,108$ virtual detectors were distributed over the phantom surface to record the photon fluence rates. Then two light sources of radius $1.0$ mm were embedded into the phantom, and filled with nanophosphors of concentrations of $5$ $\mu$g/mL and $10$ $\mu$g/mL, respectively. The centers of the two light sources were at $(2.5, 2.5, 0.0)$ and $(3.5, 0.0, 0.0)$. Two spherical subregions are nearly connected to test the spatial resolution. Using a polycapillary lens, X-ray beams are reshaped to double cone with a cone angle of $19.2$ degrees. The focal point of the X-ray is scanned at grid points in a region of interest  $\text{ROI} = \{ -2<x<5,\; -2<y<5,\; -2<z<2   \} $  along a scanning direction.  Then, the phantom was rotated with $10$ equal angles over a $360$ degrees range to acquire sufficient information for an improved stability of the reconstruction.  
The stability condition \r{U2} then holds, see Remark~\ref{rem1}.

\begin{figure}[!ht]   
  \centering
  \includegraphics[width = 2in]{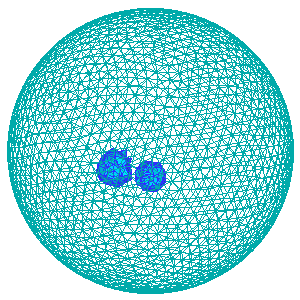} \ \ \ \ \ \ \ \ \ \     
  \includegraphics[width = 2.0in]{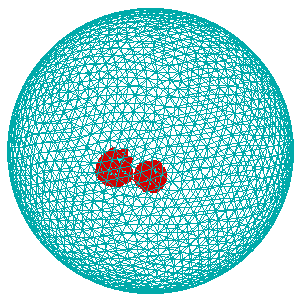}
  \caption{Modulated luminescence tomography simulation. Left: the true nanophosphor distribution in the phantom. Right: the reconstructed nanophosphor distribution. 
 The different colors here do not represent different values.  }
  \label{pic1}
\end{figure}

\begin{figure}[!ht] 
  \centering
  \includegraphics[width = 2.8in]{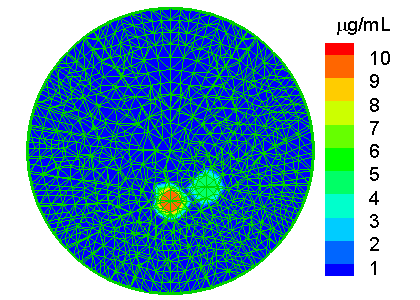} \ \    
  \includegraphics[width = 2.8in]{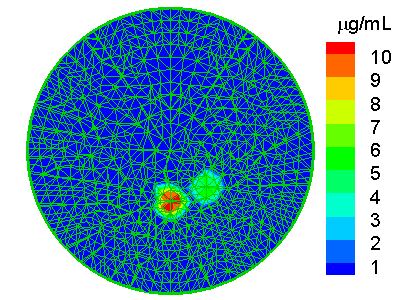}
  \caption{Modulated luminescence tomography simulation. Left: The true nanophosphor distribution on the 2D slice  $z=0$. Right:  the reconstructed nanophosphor distributionat on the 2D slice $z=0$.  }
  \label{pic2}
\end{figure}

At each scanning, the intensity of the photon luminescence on the surface of the phantom was acquired. The intensity describes the optical emission of nanophosphors with double cone excitation.  
 Poisson noise was added to the synthetic data for the simulation of measurements. After discretization,  the LSQR method \cite{Paige:1982} was used to solve solve the resulting system, in order to reconstruct the nanophosphor distribution from simulation data. The reconstructed results are in excellent agreement with the true phantom, and the average relative error of the reconstructed nanophosphor concentration was less than $5.37\%$, which was defined as 
 \[
 \text{Error} = \frac1{ \#\{i|\; \rho_i^T>\epsilon  \}}   \sum_{\rho_i^T>\epsilon} \frac{\rho_i^r-\rho_i^T }{  \rho_i^T}    ,
  \] 
where $ \rho_i^T $ and $\rho_i^r$ are the true and reconstructed nanoparticles concentrations, respectively,  and $\epsilon$ is a background noise level.  Figures~\ref{pic1}--\ref{pic2} presents the comparison between the true and reconstructed nanophosphor distribution, showing the quantification accuracy of the image reconstructions. 

\bibliographystyle{abbrv}
\bibliography{myreferences,Plamenbiomed}
\end{document}